\documentclass[12pt]{article}

\usepackage{amsmath}
\usepackage{amssymb}
\usepackage{graphics}
\usepackage{epsfig}

\begin{document}

\baselineskip 18pt

% ================================================================
% Document specific macros
% ================================================================

\newcommand{\x}{{\hat{x}}}
\newcommand{\y}{{\hat{y}}}
\newcommand{\R}{{\hat{R}}}

\newcommand{\brk}[1]{\left(#1\right)}   
\newcommand{\Brk}[1]{\left[#1\right]}   
\newcommand{\BRK}[1]{\left\{#1\right\}} 
\newcommand{\Abs}[1]{\left| #1 \right|} 
\newcommand{\pd}[2]{\frac{\partial#1}{\partial#2}}
\newcommand{\deriv}[2]{\frac{d#1}{d#2}}
\newcommand{\Bold}[1]{\boldsymbol{#1}}

\newcommand{\B}{{\cal B}}
\newcommand{\bphi}{{\varphi}}
\newcommand{\calH}{{\mathfrak H}}
\newcommand{\calR}{{\mathfrak R}}
\newcommand{\define}{\stackrel{\mathrm{def}}{=}}
\newcommand{\E}{{\mathbb E}}
\newcommand{\G}{{\boldsymbol \Gamma}}
\renewcommand{\H}{{\cal H}}
\newcommand{\Ham}{\mathfrak{H}}
\renewcommand{\O}{{\cal O}}
\renewcommand{\P}{{\cal P}}
\newcommand{\peq}{p_{\text{eq}}}
\newcommand{\s}{\Bold{s}}
\newcommand{\sig}{\mathfrak{G}}
\newcommand{\Sperp}{S_\perp}
\renewcommand{\u}{\Bold{u}}
\renewcommand{\v}{\Bold{v}}
\newcommand{\w}{\Bold{w}}

\makeatletter
\@addtoreset{equation}{section}
\makeatother
\renewcommand{\theequation}{\arabic{section}.\arabic{equation}}

% ================================================================

% ================================================================

\title{ Non-Markovian Optimal Prediction}

\author{
Alexandre J. Chorin, Ole H. Hald\\
Department of Mathematics \\
Lawrence Berkeley National Laboratory \\
1 Cyclotron Road\\
Berkeley CA 94720
\and
Raz Kupferman\\
Institute of Mathematics\\
The Hebrew University\\
Jerusalem, Israel}

\date{}
\maketitle

% ================================================================
\begin{abstract}
Optimal prediction methods compensate for a lack of resolution in the
numerical solution of complex problems through the use of prior
statistical information. We know from previous work that in the presence
of strong underresolution a good approximation needs a non-Markovian ``memory",
determined by an equation for 
the ``orthogonal", i.e., unresolved, dynamics. We present a simple approximation of
the orthogonal dynamics, which involves an ansatz and a Monte-Carlo evaluation of
autocorrelations. The analysis provides a new 
understanding of the fluctuation-dissipation formulas of statistical physics. 
An example is given. 
\end{abstract}

% ================================================================
\section{Introduction}

Many problems in science and engineering are described by nonlinear differential 
equations whose solutions are too complicated to be properly resolved.
The problem of predicting the evolution of systems that are not well
resolved has been addressed by the present authors and others in
\cite{BCC,chorin7,CKK1,CKK2,CKK3,CHK,CKL,hald,HK,K1}.  Nothing can be predicted without
some knowledge about the unresolved (``subgrid'') degrees of
freedom.  In the optimal prediction methods just cited
it is assumed that one possesses, as one
often does, prior statistical information about the system in the form
of an invariant measure; what is sought is a mean solution with
respect to this prior measure, compatible with the information
initially at hand as well as with the limitations on the computing power
one can bring to bear.

The simplest version of this idea, Markovian optimal prediction,
generates an approximating system of ordinary differential equations
and works well for a time that depends on the degree of underresolution and on the
uncertainty in the data. 
This version is optimal in the class of Markovian approximations \cite{HK}, but  
it eventually exhibits errors, because the influence of partial
initial data on the distribution of the solutions weakens in time if
the system is ergodic, and this loss of information is not captured in
full, see \cite{CHK,CKL}. To obtain an accurate approximation of a subset
of variables without solving the full problem requires the addition of a 
``memory" term, and the resulting prediction scheme becomes a generalized
Langevin equation, similar to those in irreversible
statistical mechanics \cite{evans,fick,mori,zwanzig}.

We present a general formalism for separating resolved and unresolved degrees of
freedom, analogous to the nonlinear projection formalism of Zwanzig \cite{zwanzig2} but using
the language of probability theory. We find a zero-th order approximate solution of the equation for the
orthogonal unresolved dynamics and find its statistics by Monte-Carlo integration; we use the results to
construct a prediction scheme with memory.
We apply the scheme to a simple model problem. In the conclusion we indicate how the construction is
generalized to more complicated problems, and the new perspectives it opens for prediction as well as for irreversible
statistical mechanics.

% ================================================================
\section{Linear and nonlinear projections of dynamical systems}

Consider a problem of the form
\begin{equation}
\deriv{\bphi}{t} = R (\bphi),\ \ \bphi(0)=x,
\label{eq:system}
\end{equation}
where $R$ and $\bphi$ are $n$-dimensional vectors ($n$ may be
infinite), with components $R_i$ and $\bphi_i$; $t$ is time.  When
$n$ is finite (\ref{eq:system}) is a system of ordinary differential
equations. Our goal is to calculate the average values of $m$ components of $\bphi$, $m<n,$  
without calculating all the components; the average is over all the
values that the missing, unresolved, components may assume; our prior information allows us
to make statistical statements 
about these missing components. 

We denote the phase space (the vector space in which $\bphi$ resides) by
$\G$; in classical statistical physics this phase space is the $n=6
\ell$ dimensional space of coordinates and momenta $(q_i,p_i)$, where
$\ell$ is the number of particles; the $q_i,p_i$ at time $t$ are then entries of the vector $\bphi$.  A solution of equation
(\ref{eq:system}) is defined when an initial value $\bphi(t=0) = x$ is
given; to each initial condition $x$ corresponds a trajectory,
$\bphi(t)=\bphi(x,t)$; the initial value $x$ is emphasized by this notation in view of its key role in what follows.  

A phase variable $u$ is a function on $\G$; $u$ may be a vector, whose components
are labeled as $u_i$.
A phase variable varies when its argument varies in time, so that a
phase variable whose value at $t=0$ was $u(x)$ acquires at time $t$
the value $u(\bphi(x,t))$. It is useful to examine the evolution of $u$
in a more abstract setting:
Introduce an evolution operator $S^t$ for phase variables by the relation
\begin{equation}
(S^t u)(x)= u(\bphi(x,t)).
\label{eq:pullback}
\end{equation}
Differentiation of (\ref{eq:pullback}) with respect to time 
yields
\begin{equation}
\pd{}{t} (S^t u)(x) = 
\sum_i R_i(\bphi(x,t)) \pd{u}{\bphi_i}(\bphi(x,t)) =
L S^t u(x),
\label{liouville}
\end{equation}
where $L$, the Liouvillian, is the linear differential operator
$L=\sum_i R_i(\bphi) \pd{}{\bphi_i}$.  Thus the phase variable $S^tu$ can be
calculated in either of two ways: (i) for each $x$ integrate to time $t$ the
equations of motion $\deriv{}{t}\bphi(t) =
R(\bphi(t))$ with initial conditions $\bphi(0)=x$
and evaluate the phase variable $u$ at the point
$\bphi(x,t)$; or (ii) solve the equation
\begin{equation}
\left\{
\begin{split}
& \pd{}{t} S^t u = L S^t u      \\
& S^0 u = u(x).
\end{split} \right.
\nonumber
\end{equation}
It is convenient to write $S^t=e^{tL}$;
we do not inquire here as to the conditions under which this symbolic notation
can be taken literally. The significant thing about equation (\ref{liouville}) is that it
is linear.    

One can check from the definitions 
that $e^{tL}(uv) = (e^{tL}u) (e^{tL}v)$, $e^{tL}f(u) = f(e^{tL}u)$ for any function $f$, and
$e^{tL}L=Le^{tL}$. In this notation, the symbol $u$ standing alone refers to the data $u(x)$ at $t=0$; the time dependence
is described by the exponential. 

Suppose that the initial data $x$ are drawn from a probability
distribution $\mu^0$; each initial datum gives rise to a solution of
equation (\ref{eq:system}) and the measure $\mu^0$ evolves into a
measure $\mu^t$ at time $t$. The evolution of $\mu^t$ is defined by
the conditions
\begin{equation}
\int_{\G} u(\bphi(x,t)) \,\mu^0(dx) = \int_{\G} u(x) \,\mu^t(dx)    
\nonumber
\end{equation}
for all sufficiently smooth phase variables $u$.  We assume that the
measure $\mu^0$ is invariant under the flow (\ref{eq:system}): $\mu^t
= \mu^0=\mu$. Many systems have invariant measures, in particular, any Hamiltonian
system with Hamiltonian $H$ leaves invariant the canonical measure with density $Z^{-1}e^{-H(x)/T}$, where $Z$ is a normalization constant and $T$ is the variance of the samples, which in physics is the temperature.
Given a phase variable $u$, we denote by $\E[u]$ the
expected value of $u$ with respect to the invariant measure $\mu$,
\begin{equation}
\E[u] =\int_\G u(x)\,\mu(dx).
\nonumber
\end{equation}
We endow the space of phase variables with the inner product $(u,v)=\E[uv]$, which 
makes them elements of the Hilbert space $L_2[\G,\mu]$ ($L_2$ for brevity). 
The representations (\ref{eq:system}) and (\ref{liouville}) of the dynamics are equivalent; in particular one can retrieve (\ref{eq:system}) from (\ref{liouville})
by setting $u=x_i$ and therefore $e^{tL}u=\bphi_i(x,t)$ (the i-th coordinate function).
In $L_2$ the Liouvillian operator $L$ is antisymmetric.

We wish to calculate the means of a small number of variables $\bphi_i(x,t), i=1,\dots,m, m<n$,
averaged over all the values the other variables may take initially when they are 
drawn from the invariant distribution, 
without calculating any of these other variables.
We write $\hat\bphi=(\bphi_1,\dots,\bphi_m)$
for the vector whose entries are the variables we actually calculate. 
To find equations for $\hat\bphi$ we need projections of $L_2$ on subspaces of
functions of ${\hat x}=(x_1,\dots,x_m)$ or of $\hat\bphi(x,t)=(\bphi_1,\dots,\bphi_m)$; these projections are not unique and
we shall consider two different ones:

(i) 
Consider the 
conditional expectation $\E[v|u]$, where both $u$ and $v$ are phase
variables; it satisfies:
\begin{enumerate}
\item
$\E[v|u]$ is a function of $u$;
\item
$\E[v|u]$ is linear in $v$:
\begin{equation}
\E[\alpha v_1 + \beta v_2|u] = \alpha\,\E[v_1|u] + \beta\,\E[v_2|u].
\nonumber
\end{equation}
\item
$\E[v|u]$ is the best approximation of $v$ by a function of $u$:
\begin{equation}
\E[|v-\E[v|u]|^2] \le \E[|v-f(u)|^2]
\label{best}
\end{equation}
for all functions $f$.  
\end{enumerate}
See Chung \cite{chung}.

$\E[v|u]$ is the orthogonal
projection of $v$ on the space of functions of $u$, and 
we can write
$Pv=\E[v|u]$. This is the ``nonlinear projection"
used in \cite{zwanzig2} with a different interpretation, as well as in \cite{CKK1,CKK2,CHK,CHK}.

ii)For a function $w=w(x)$ in $L_2$, define
\begin{equation}
P'w=\Sigma_1^m\frac{(w,x_i)}{(x_i,x_i)}x_i,
\label{linearproj}
\end{equation}
this is the ``linear projection" used in most of irreversible statistical mechanics (see e.g. (\cite{fick,forster,evans})
.

Clearly the projection based on conditional expectations puts more information into
the subspace of functions of the smaller set of variables and is thus preferable for our purposes,
but we shall need both projections for technical reasons. 
There is a large set of projections intermediate between these two, obtained by spanning the
space of functions of $\x$ by additional elements; we shall not discuss these intermediate projections here.

We now follow the Mori-Zwanzig procedure (\cite{evans,forster,mori, zwanzig}) and
split the time derivative of $e^{tL}u$ into its
projection on the functions of $\hat\bphi(t)$ plus a complement:
\begin{equation}
\pd{}{t} e^{tL}u = e^{tL} Lu =  e^{tL}PLu+e^{tL}QLu, 
\label{eq:split}
\end{equation}
where $P$ is either of the two projection just described and $Q=I-P$. 
If $u$ is a scalar function of $\bphi$, it is projected in $L_2$ on a function of
$\hat\bphi$; if $u$ is a vector function of $\bphi$, each of its components is projected. 
The first term in (\ref{eq:split}) can be evaluated in terms of $u$ and we 
denote it by $\calR(u)$: 
\begin{equation}
e^{tL}\calR(u) = \calR(e^{tL}u).
\nonumber
\end{equation}

To understand the second term, consider an evolution operator,
$e^{tQL}$, defined for a phase variable $v$ by
the equation:
\begin{equation}
\pd{}{t} e^{tQL}v  =QLe^{tQL}v= Le^{tQL}v - PLe^{tQL}v. 
\end{equation}
with $Pv(t=0)=0$. This is the orthogonal dynamics equation, which should be read as 
follows: Let $w=v$ at $t=0$; let $w(t)$ be a time-dependent phase variable that satisfies
\begin{equation}
\pd{}{t}w(t)=Lw(t)-PLw(t);
\end{equation}
with initial datum $w(t=0)=v$, 
and set $w(t)=e^{tQL}v$, defining the orthogonal evolution operator $e^{tQL}$. 
Note that 
if $v$ is orthogonal to the range of $P$ then so is $e^{tQL}v$ for 
all times $t$. The operator $e^{tQL}$  
is the solution operator of the orthogonal dynamics.

The evolution operators $e^{tL}$ and $e^{tQL}$ satisfy the Dyson formula
\cite{evans}:
\begin{equation}
e^{t L} = e^{tQL} + \int_0^t e^{(t-s) L} P L e^{s QL} \,ds,
\label{dyson}
\end{equation}
which can be checked by differentiation. 
With the help of the Dyson
formula the second term on the right hand side of (\ref{eq:split}) can
be written as:
\begin{equation}
e^{tQL}v+\int_0^te^{(t-s)L}PLe^{sQL}vds,
\nonumber
\end{equation}
where 
\begin{equation}
v=Lu-PLu=QLu.
\nonumber
\end{equation}

Putting all the terms together, we obtain the generalized Langevin
equation
\begin{equation}
\pd{}{t}e^{tL}u  =  \calR(e^{tL}u)+
\int_0^t e^{(t-s)L}PLe^{sQL}vds + e^{tQL}v.
\label{eq:langevin}
\end{equation}
This is  an identity between phase variables, the
starting point for our
approximations.
If one projects this equation on the space spanned by the initial data
for $u$, the last term drops out and the other terms acquire a projection operator
as prefactor. 

The various terms in equation (\ref{eq:langevin}) have conventional
interpretations.  The first term on the right-hand side is a function
only of $e^{tL}u$ and represents the self-interaction of the
components of $e^{tL}u$; it is the Markovian contribution to
$\pd{}{t} e^{tL}u$.  The second term depends on $u$ through the
values of $e^{tL}u$ at all times $s$ between $0$ and $t$, and
embodies a non-Markovian memory.  Finally, the third term has no
component in the range of $P$ and if we know nothing outside this
range it can be viewed as random, with statistics determined by the
initial data.

In the special case of a linear projection $P'$ the second term in equation (\ref{eq:langevin})
acquires a particularly simple form, well-known in statistical physics where it is the only
form in general use.  
By definition, see above,
\begin{equation}
P'Le^{sQL}QLu=P'LQe^{sQL}QLu= \sum_1^m(LQe^{sQL}QLu,h_j(\x))h_j(\x),
\end{equation}
where $h_j(\x)=\x_j/||\x_j||$, $||\x_j||^2=E[\x_j^2]$, $1\le j\le m$, are normalized coordinate functions.
A short manipulation converts this sum into
\begin{equation}
-\sum K_j(s)h_j(\x),
\end{equation}
where $K_j(s)=(e^{sQL}QLu,QLh_j)$, i.e., $K_j$ is the  
correlation of the solution of the orthogonal dynamics equation
that starts from $v=QLu$ with $QLh_j$. Substitution into the integral  
produces the term
\begin{equation}
-\int_0^t\sum K_j(t-s)e^{sL} h_j(\x)ds.
\label{kernels}
\end{equation} 

% ===================================================================

\section{Implementation of the projections in a small dimensional space and 
first-order optimal prediction}

The generalized Langevin equation (\ref{eq:langevin}) expresses the
rate of change of the phase variable $u$ as a sum of terms that 
depend on $u$ and 
on the orthogonal dynamics. These expressions cannot be used directly for approximation
when the equations are nonlinear. Indeed, the evaluation of a term
such as $e^{tL}E[Lu|u]=E[Le^{tL}u|e^{tL}u]$ requires the evaluation of
$e^{tL}u$, i.e., requires a solution of the full set of $n$
equations. In terms of $\bphi$, this term becomes

\begin{equation}
\E[\R(\y)|\y]_{\y=\hat\bphi(x,t)},
\end{equation}
where the last expression is the expected value of $\R(\y)$ given the value $\hat\bphi(x,t)$  
of $\y$. 
(
we first find the best approximation of $\R(y)$ by a function of $\y$, then
substitute the value $\y=\hat\bphi(x,t)$ into the result).
The vector $x$ has components $\tilde{x}$ not included in
$\x$ and not known; $x=(\x,\tilde{x})$.  To obtain something that
can be evaluated we move the expectation into the condition, i.e., we
approximate $\E[\R|\y]_{\y=\bphi(\x,\tilde{x},t)]}$ by
$\E[\R|\y]_{\y=\E[\bphi(x,t)|\x]}$, which is
a function $\y=\y(\x)$ of $\x$.

Assume for a moment that the second and third terms in the Langevin equation (\ref{eq:langevin}) are
negligible. Equation  (\ref{eq:langevin}) becomes:

\begin{equation}
\frac{d}{dt}\hat\bphi(x,t)
=\E[\R(y)|\y]_{\y=\hat\bphi(x,t)}.
\label{outer}
\end{equation}
The approximation we just made, replacing $\hat\bphi(x,t)$ by its projection, 
gives:
\begin{equation}
\frac{d}{dt}\hat\bphi(\x,t)=\E[\hat\R(\y)|\y]_{\y=\y(\x)};
\end{equation}
projecting this equation on the space of functions of $\x$, we find,
\begin{equation}
\frac{d}{dt}\y=\E[\R(\y)|\y]_{\y=\y(\x)}.
\label{op}
\end{equation}
The initial data $\y(0)=\x$ make $\y(t)$ a function of $\x$.  
These are the first-order optimal prediction equations. 
First-order optimal prediction neglects however the
effect of the fluctuations due to the variation of the components $\tilde{x}$ of $x$, as we already know from
\cite{CKL}; all that is left the uncomputed components is their conditional expectations.

As an example, consider the model equations that arise from the Hamiltonian
\begin{equation}
H=\frac{1}{2}(x_1^2 + x_2^2 + x_3^2+ x_4^2 + x_1^2 x_3^2),
\end{equation}
($x_1,x_3$ can be viewed as position variables $q_1,q_2$ and $x_2,x_4$ can be viewed as momenta $p_1,p_2$). The resulting
equations of motion are:
\begin{eqnarray*}
\frac{d}{dt}\bphi_1&=&-\bphi_2-\bphi_2\bphi_4^2,\\
\frac{d}{dt}\bphi_2&=&\bphi_1,\\
\frac{d}{dt}\bphi_3&=&-\bphi_4-\bphi_4\bphi_2^2,\\
\frac{d}{dt}\bphi_4&=&\bphi_3
\end{eqnarray*}
with $\bphi(0)=x$. 
The probability density of $x$ is $Z^{-1}e^{-H(x)/T}$, and that of $\bphi$ is $Z^{-1}e^{-H(\bphi)/T}$.
For simplicity set $T=1$. 
Keeping $m=2$ out of the $n=4$ equations, and omitting all but the first term in the Langevin
equation (\ref{eq:langevin}), we find:
\begin{eqnarray*}
\frac{d}{dt}\bphi_1&=&-\bphi_2-\frac{\bphi_2}{1+\bphi_2^2},\\
\frac{d}{dt}\bphi_2&=&\bphi_1
\label{toosimple}
\end{eqnarray*}
as can be deduced from the definition of conditional expectation:
\begin{equation}
\frac{\int_{-\infty}^{+\infty}x_2x_4^2 e^{-H} dx_3 dx_4} 
{\int_{-\infty}^{+\infty} e^{-H}dx_3dx_4}=\frac{x_2}{1+x_2^2}. 
\end{equation}
These are the first-order optimal prediction 
equations as presented in \cite{CKK1,CKK2,CKK3,CHK,CKL}.

\section{Zero-th order solution of the Langevin equation}

Now consider the full Langevin equation (\ref{eq:langevin}). 
The evaluation of the memory term requires a solution of the
orthogonal dynamics equation, which we shall present in detail
elsewhere. In the present paper we shall be content with less,
yet the result is instructive.

With the appropriate substitutions, the Dyson formula (\ref{dyson})
yields
the identity
\begin{equation}
e^{tQL}=e^{tL}-\int_0^t e^{(t-s)L}PLe^{sQL}ds,
\label{inversedyson}
\end{equation}
which can be in principle be the starting point of a perturbative evaluation
of $e^{tQL}$. 
The zero-th order approximation of this equation is:
\begin{equation}
e^{tQL}\cong e^{tL},
\end{equation}
i.e., one replaces the true flow in the orthogonal complement of the range of
$P$ (or $P'$) by the real flow induced by $L$ ( which of course has a component in
the orthogonal complement). 
This approximation can be sufficient if the correlations in (\ref{kernels})
do not decay too slowly; to see this, consider the second term in equation (\ref{eq:langevin}):

\begin{equation}
\int_0^te^{(t-s)L}PLe^{sQL}QLuds=\int_0^te^{(t-s)L}PLQe^{sQL}QLuds.
\end{equation}
Adding and subtracting equal quantities, we find:
\begin{equation}
PLe^{sQL}QLu=PLQe^{sL}u+PLQ(e^{sQL}-e^{sL})QLu;
\end{equation}
a Taylor series yields:
\begin{equation}
e^{sQL}-e^{sL}=I+sQL+\dots-I-sL-\dots=-sPL+O(s^2),
\end{equation}
and therefore, using $QP=0$, we find:
\begin{eqnarray*}
PLe^{sQL}QLu&=&PLQe^{sL}QLu+O(s^2),\\
\int_0^te^{(t-s)L}PLe^{sL}QLuds&=&\int_0^te^{(t-s)L}PLe^{sQL}QLuds+O(t^3);
\end{eqnarray*}
if the correlations $(e^{sQL}QLu,QLh_j)$, $(e^{sL}QLu,QLh_j)$, are significant only over
short times, we have an acceptable approximation. 

At $t=0$ the orthogonal dynamics start from $v(0)=R(x)-\E[R(x)|\x]$;
at later time, in our approximation,
$v(t)=e^{tL}v(0)=R(\bphi(x,t))-\E[R(y)|\y]_{\y=\hat\bphi(x,t)}$.
The correlations in (\ref{kernels}) can be calculated by Monte-Carlo, and  
provides a better approximation than first-order optimal prediction. This is still not
a very good approximation, and the reason, as we shall show elsewhere, is that the linear
projection $P'$ used to derive (\ref{kernels}) projects onto a set of functions that fails to span 
the whole space of functions of $\x$; this can be remedied by devising more elaborate approximations to the
conditional expectation $P$. 
Here we propose a heuristic alternative: 
We interpret $\E[R(y)|\y]_{\y=\hat\bphi(x,t)}$ in $v(t)$ as the mean of the right-hand-side of equations (\ref{eq:system}) and $R(\bphi(x,t))$
as the fluctuation around that mean. This mean varies less than the fluctuations as $x$ spans $\G$,
and we make it vary even less by anchoring it to the initial data for the specific problem
we wish to solve, i.e., first approximate 
$\E[R(y)|\y]_{\y=\hat\bphi(x,t)}$ by $\E[R(y)|\y]_{\y=\E[\hat\bphi(x,t)|\x]}$ 
and then further fix $\x$ at the specific initial value for which we want to solve the $m-$equation approximation of the  system (\ref{eq:system}).
Let $w(t)$ be the function $w(t)=\E[R(y)|\y]_{\y=\E[\hat\bphi(x,t)|\x]}$
obtained in this way for
a specific $\x$ (in the example, $w(t)=-\bphi_2-\bphi_2/(1+\bphi_2^2)$, where $\bphi_2$ is the specific
solution we are seeking). Note that $v(t)=R(\bphi(x,t))-w(t)$ is no longer a stationary stochastic process. 
The heuristic rationale for this ``freezing" of $w$ is as follows: Write as before $x=(\x,\tilde{x})$; 
when, in the averaging that determines the kernels $K_j$, we fix $\x$ and vary $\tilde{x}$ alone, thus finding the contribution of the subspace orthogonal to $\x$, the ``frozen" and the true variations of $v(t)$ are the same;
however, if we keep $\tilde{x}$ fixed and vary $\x$ for the evaluation of $R(\bphi)$ we relate the fluctuations to the
specific solution we are seeking more thoroughly than is possible otherwise when the projection is $P'$.   
With these approximations, 
the kernels $K_j$ in (\ref{kernels}), $K_j(t-s)=\E[v(t)v(s)]$, become
\begin{equation}
K_j(t-s)=K_j^0(t-s)-w_j(t)w_j(s).
\label{haha}
\end{equation}
where $K_j^0=\E[R_j(\bphi(x,t)R_j(\bphi(x,s)]$ 
can be evaluated by Monte-Carlo: Sample initial data 
from the invariant distribution over and over, solve the resulting system of equations,
whose solutions are statistically stationary in time, and average. Note that
this calculation does not depend on the particular known initial values $\x$,  and thus can be done
once and for all for a given equation and a given level of truncation. The time evolution of the $v_j(t)$ need be accurate only as long as $v_j(t)$ and $v_j(0)$ are correlated. 

Substitute this into the Langevin equation: there is no impediment to using the true
$e^{tQL}$ in the third term; then average with respect to $x$; 
note that the 
$K_j$ are functions of time only  and commute with a projection on the initial data,
while the average of the third term is 0. 
We find the following equation:
\begin{equation}
\frac{d}{dt}\y_j=w(t)-\int_0^t\sum_1^m(K_j^0(t-s)-w_j(t)w_j(s))\y_j(s)ds,
\label{oops}
\end{equation}
where $j=1,\dots,m$ and 
$\y(\x,0)=\x$.

Equation (\ref{oops})  
states that the solution of the problem depends not only on the mean of the $R$ but also on the autocorrelation of the fluctuations. The solution is thus 
``renormalized" ( see e.g \cite{hou}). The sources of error in (\ref{oops}) are: An inconsistent
use of projections, a simplistic solution of the orthogonal dynamics, and possibly an
inaccurate Monte-Carlo evaluation of the kernels $K^0$. All can be remedied. 

In the example, (\ref{oops}) reduces to:
\begin{eqnarray}
\frac{d}{dt}y_1&=&w(t)-\int_0^t(K^0(t-s)-w(t)w(s))y_1(s)ds,\\
\frac{d}{dt}y_2&=&-y_1,\hfill
\label{example2}
\end{eqnarray}
The kernel $K^0$ is $K^0=\E[R_1(\bphi(x,t))R_1(\bphi(x,s))], R_1(\bphi)=
-\bphi_2-\bphi_2\bphi_4^2$, $w(t)=-y_2-y_2/(1+y_2^2)$, and the expected value in $K^0$ is evaluated over all choices of $x$ in $n=4$ 
dimensions drawn from the invariant distribution with density $Z^{-1}e^{-H/T}$.

In Figure 1 we display some numerical results for $y_1(t)$ in the example, with  data $x_1=1,x_2=0$.
We show the ``truth" found by averaging many solutions of the full $4-$dimensional system, 
the Galerkin approximation that sets all unknown functions to zero, the first-order optimal
prediction, and the solution of equation (\ref{example2}). The solutions  are shown in a rather favorable
case ( they look less striking if one exchanges for example the values of $x_1,x_2$).

\section{Conclusions}

At first sight, the results in Figure 1 are interesting but not overwhelming: The cost of
evaluating $K^0$ is comparable to the cost of evaluating the "truth" by Monte-Carlo, and
the gain is not obvious. One may note however, that once $K^0$ has been evaluated, 
the cost of 
rerunning the calculation with any other initial data $\x$ is negligible. 

This is the significant fact:  $K^0$ is evaluated ``at equilibrium", i.e., with all
components of $x$, the initial data, sampled from the invariant measure; $K^0$
does not depend on any specific initial value. As we shall show elsewhere, the analogous
statement is true for an accurate solution of the orthogonal dynamics. Once the 
heavy work of determining memory functions has been done, the solution of a specific
problem is plain sailing. At equilibrium, one can bring to bear the panoply of scaling methods
and equilibrium statistical mechanics. One can say that the Mori-Zwanzig formalism 
makes possible the use of ``universal" (non problem specific) results to solve specific problems.
In some of the applications we have in mind, the ``large" ($n$-dimensional) problems are partial differential
equations, and then imperfections in the evaluation of memory terms are immaterial as long as  
the rate of convergence of finite-dimensional approximations is enhanced ( see \cite{BCC,CKK2,CKL}).

It is taken for granted in the physics literature that the memory kernels are autocorrelations
of the ``noise" (i.e., the orthogonal dynamics). This is true also in the example we have presented here. However, it should be obvious
from the discussion that this is an artifact of the use of $P'$, the linear projection,
and that the full truth is more complicated and interesting.

Furthermore, in the physics literature one usually deals with memory by separating ``fast" and "slow"
variables and assuming that the orthogonal dynamics generated by ``fast" variables generate ``noise" with delta-function memory ; note in contrast  that in our problem the unresolved and the resolved
variables have exactly the same time scale. 

Finally, the heuristic arguments of the present paper will be replaced in general by a systematic approximation of
the Langevin equation, including a systematic evaluation of the orthogonal dynamics, as we will explain in
future publications. 

\begin{figure}[htbp] 
\epsfig{file=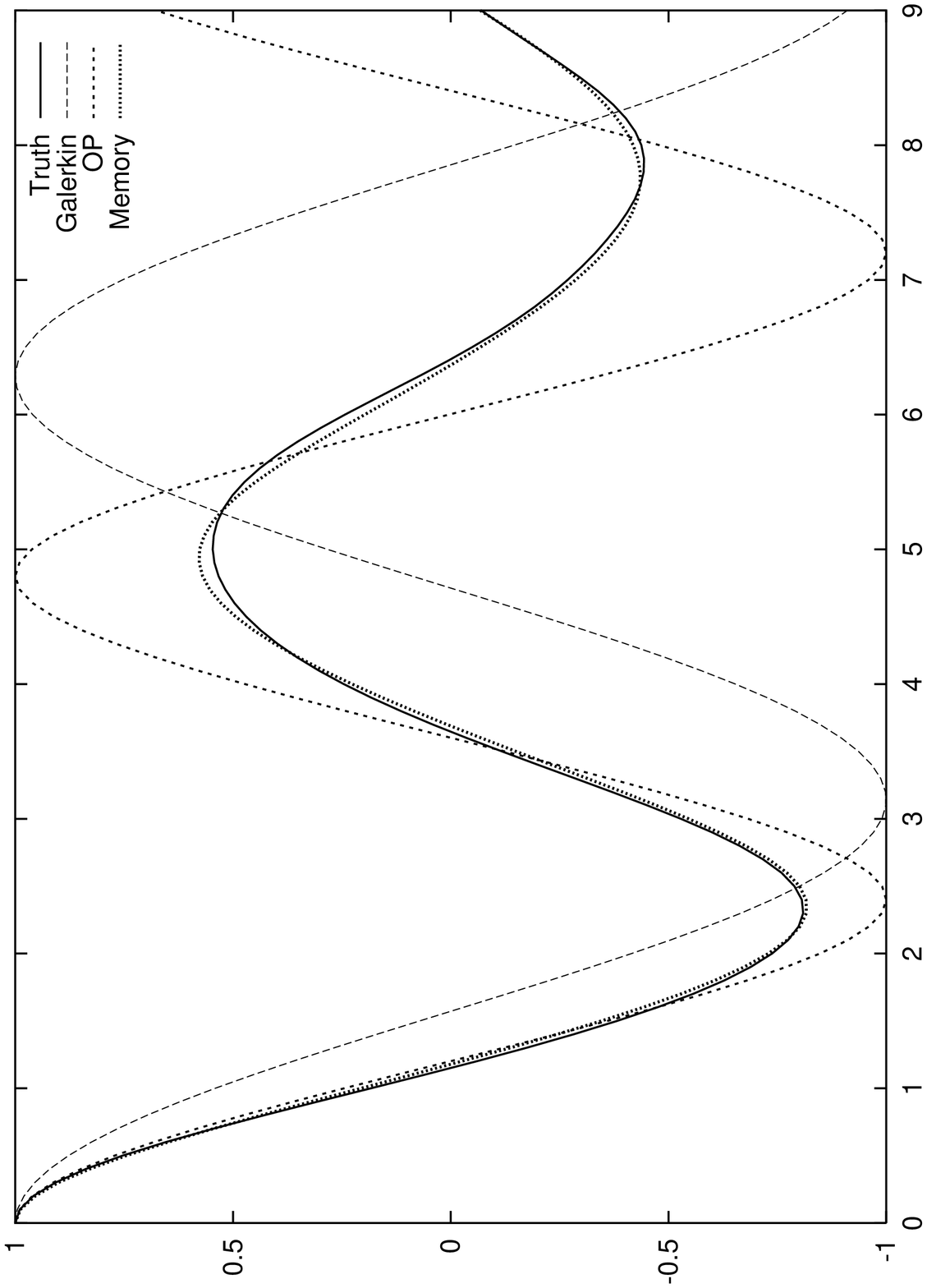,width=4.5in}
\centerline{Figure 1}
\end{figure}

\section{Acknowledgements} We would like to thank Prof. G.I. Barenblatt, Dr. E. Chorin, Mr. E. Ingerman, Dr. A. Kast, Mr. K. Lin, for helpful discussions and comments, and Mr. P. Okunev for help in programming. 
This work was supported in part by by the Applied Mathematical Sciences subprogram of the Office of
Energy Research of the US Department of Energy under Contract DE-AC03-76-SF00098,
and in part by the National Science Foundation under Grant DMS98-14631. 
R.K. was supported by the Israel Science Foundation founded by the Israel Academy of Sciences and Humanities.

\end{document}